\documentclass[12pt]{article}
\usepackage{amsmath}
\usepackage{amssymb}
\usepackage{amsfonts}
\usepackage{amsthm}
\usepackage{latexsym}
\usepackage{epsfig}
\usepackage{array}
\usepackage{boxedminipage}
\makeatletter
\renewcommand{\section}{\@startsection
{section} {1} {0mm} {-\baselineskip} {0.5\baselineskip}
{\large\bf}}
\renewcommand{\subsection}{\@startsection
{subsection} {2} {0mm} {-\baselineskip} {0.5\baselineskip}
{\normalsize\bf}} \makeatother

\newtheorem{defin}{Definition}
\newtheorem{teorema}{Theorem}

\newcommand{\ke}{\mathrm{ker} \hspace{0.5mm}}

\begin{document}
%\title{}
%\author{Ricardo Riaza
%\thanks{The author is with the Departamento de %Depto.\ 
%Matem\'{a}tica
% Aplicada a las Tecnolog\'{\i}as de la Informaci\'{o}n,
%Escuela T\'{e}cnica Superior de Ingenieros de 
%ETSI Telecomunicaci\'{o}n,
%Universidad Polit\'{e}cnica de Madrid  %\\ Ciudad Universitaria s/n 
%- 28040 Madrid, Spain. E-mail: {\sl rrr@mat.upm.es}.
%Supported by Research Project MTM2007-62064
%of Mi\-nis\-terio de Educaci\'{o}n y Ciencia, Spain.}}
%\maketitle

%\markboth{IEEE Transactions on Circuits and Systems --} 
%~Vol.~6, No.~1, January~2007}%
%{Riaza:}

%\mbox{}

\begin{center}
{\Large\bf  %Hybrid memristors and the frontiers of %, f
First order devices, hybrid memristors, and \vspace{2mm} \\
the frontiers of 
%classical 
nonlinear circuit theory\footnote{A revised version
of this manuscript has been published as R. Riaza,
{\em First Order Mem-Circuits: Modeling,
Nonlinear Oscillations and Bifurcations} in the
IEEE Transactions on Circuits and Systems I, Vol. 60, no. 6, June
2013, pp. 1570-1583. DOI: 
{\sf http://dx.doi.org/10.1109/TCSI.2012.2221174}.}\footnote{Supported by 
%{\Large\bf  Generalized circuit theory: \vspace{3mm} \\
%Hybrid memristors and first order devices\footnote{Supported by 
%{\Large\bf  Generalized circuit theory: %\vspace{3mm} \\
%First order devices\footnote{Supported by 
%{\Large\bf  Beyond memdevices:
%Generalized circuit theory\footnote{Supported by 
Research Projects MTM2007-62064
of Mi\-nis\-terio de Educaci\'{o}n y Ciencia
and MTM2010-15102
of Mi\-nis\-terio de Ciencia e Innovaci\'{o}n, Spain.}}

\ \vspace{1mm} 
\\

\ \\

%{\Large\bf Generalized circuit theory: \vspace{1mm}\\
%higher order devices and mem-systems\footnote{Supported by 
%Research Project MTM2007-62064
%of Mi\-nis\-terio de Educaci\'{o}n y Ciencia, Spain.}}

%\ \vspace{-2mm} \\

%\ \\
%\setcounter{footnote}{6}

{\large\sc Ricardo Riaza}\\ %\footnote{Corresponding author}}\\
\ \vspace{-3mm} \\
Departamento de Matem\'{a}tica Aplicada a las Tecnolog\'{\i}as de la
Informaci\'{o}n \\
Escuela T\'{e}cnica Superior de Ingenieros de Telecomunicaci\'{o}n \\
Universidad Polit\'{e}cnica de Madrid  %\\ Ciudad Universitaria s/n 
- 28040 Madrid, Spain \\
{\sl ricardo.riaza@upm.es} \\

%\vspace{0.5cm}

%[DRAFT]

%\

%\

\end{center}

\

\begin{abstract}
Several devices exhibiting memory effects have shown up in 
nonlinear circuit theory in recent years. Among others, these circuit elements 
include Chua's memristors, as well as memcapacitors and
meminductors. These and other related devices 
seem to be beyond the, say, classical scope of circuit theory,
which is formulated in terms of 
resistors, capacitors, inductors, and voltage and current
sources.
We explore in this paper the potential extent of nonlinear circuit theory
by classifying such mem-devices in terms of the variables 
involved in their constitutive relations and the notions of the differential-
and the state-order of a device. Within this framework, 
the frontier of first order circuit theory is
defined by so-called {\em hybrid memristors}, which are proposed here
%in order 
to accommodate a characteristic relating all four 
fundamental circuit variables.
Devices with differential order two 
and mem-systems are discussed in less detail.
We allow for fully nonlinear characteristics in all circuit elements,
arriving at
a rather exhaustive taxonomy of $C^1$-devices.
Additionally, we extend the
notion of a topologically degenerate configuration to circuits with
memcapacitors, meminductors and all types of memristors, 
and characterize the
differential-algebraic index of nodal models of such circuits.
%, extending previous results.
\end{abstract}

%\vspace{5mm}

\

\noindent {\bf Keywords:} nonlinear circuit, memristor, 
%nodal analysis, branch-oriented model,
%state model, 
memcapacitor, meminductor, nodal analysis,
differential-algebraic equation, index.
%semistate model, 
%matrix pencil, 
%index, tree.
%circuit model, tableau approach,
%modified node analysis, augmented node analysis,
%graph topology, 
%singular ODE, bifurcation, 
%matrix pencil. 
%state-space equation.

\

%\vspace{5mm}

%\

\noindent {\bf AMS subject classification:}
%05C50, %Graphs and matrices
%15A22, %Matrix pencils
34A09, %Implicit equations, differential-algebraic equations [See also 65L80]
%%34D20, %ODEs: stab. theory, Lyap. stab.
%%%37C75, %Smooth dyn. sys., Stab. theory
%%%34C23, %Bifurcation (ODEs)
%37G10, %Bifurcation of singular points;
94C05. %Analytic circuit theory
%94C12 Fault detection; testing
%94C15. %Applications of graph theory

%\thispagestyle{empty}

%\newpage

%\thispagestyle{empty}

%\ \mbox{}

%\newpage

%\setcounter{page}{1}

%\vspace{0.5cm}

\

\newpage
\section{Introduction}
\label{sec-intr}

Broadly speaking, nonlinear circuit theory is concerned with the study of
%(implicitly defined) 
constrained ordinary differential equations involving time and 
four $m$-dimensional variables 
$q,$ $\varphi,$ $i,$ $v$ (standing for charge, flux, current and voltage, respectively),
with the following restrictions: \label{page-abc}
\begin{itemize}
\item[(a)] the $2m$ differential relations $q'=i$, $\varphi'=v$ always hold; 
\item[(b)] the %voltage and current 
vectors $v$, $i$, satisfy a total amount of $m$ 
linearly independent relations $Bv=0$, $Di=0$ coming from Kirchhoff laws;
\item[(c)] the characteristics of devices define $m$ additional relations among
the circuit variables.
\end{itemize}
%The prime $'$ denotes time derivation.  Careful, use it later in
%another sense.
%we are assuming  an autonomous setting for simplicity in this
%introductory discussion (in particular, independent sources are
%not yet considered). 
This means that nonlinear circuit models can be generally written in the form
\begin{subequations} \label{basic}
\begin{eqnarray}
q' & = & i\\
\varphi' & = & v \\
0 & = & f(q, \varphi, i, v, t),
\end{eqnarray}
\end{subequations}
where $f$ captures Kirchhoff laws but also the constitutive relations of all
circuit devices.
Aside from the different circuit topologies, reflected in the 
form of the loop and cutset matrices
$B$, $D$ arising in item (b), the differences between circuit families
%must either 
come from %, or from 
the devices' characteristics 
referred to in (c). %We will focus our attention on (c).

%%%%%%%%%% ``Superfluous'' variables NOT HERE

In classical circuit theory,
%'' (or: Before the appearance of the memristor) 
the devices' characteristics
just involve two out of the four variables mentioned above in one of three ways,
relating current and voltage in resistors and controlled sources, 
charge and voltage in capacitors,
and flux and current in inductors. In general, these relations may be
nonlinear and/or involve time explicitly; specifically,
the current or the voltage in independent
sources is an explicit function
of time. This provides a setting which models a great variety
of circuits arising in electrical and electronic engineering; in particular,
many semiconductor devices, including transistors, may be accommodated
in this framework by means of equivalent circuits including controlled sources.
%specifically, 
%{\em resistors} are governed by a relation of the form $h_r(v_r,i_r)=0$,
%the characteristic of {\em capacitors} reads as $h_c(q_c,v_c)=0$, 
%and for {\em inductors} we have $h_l(\varphi_l,i_l)=0$. We use
%subscripts in order to distinguish the variables associated
%with different devices. *Sources: t; controlled relate i, v. Mention
%equivalent circuits...

%\newpage

\ %\vspace{3mm}

\noindent {\bf Memristors.} 
The {\em memory-resistor} or {\em memristor} is changing this picture
substantially. This device, whose existence was predicted by Leon Chua
in 1971 \cite{chuamemristor71}
and which actually appeared at the nanometer scale in 2008 
\cite{memristor2008}, is defined by a nonlinear charge-flux relation.
The potential applications of this device in the design of non-volatile
memories, signal processing, adaptive and learning systems, reconfigurable nanoelectronics, 
%field programmable gate arrays, 
etc.,
%might make it a revolutionary discovery (nmg).
might make the memristor and related devices play a very significant 
%%%%%key
role in electronics in the near future,
specially at the nanometer scale. A lot of research is focused on
this topic; cf.\
\cite{benderli09,  chen, chuamemristor08, itohchua09, kavehei10, messias10,
muthus10, muthuskokate,  pershin, pershin09, pershin10a, pershin10b,
pershin10c, 
memactive, %membranch,
fully, memristors, 
snider08, wang09, wu, yang08}.

Not only regarding applications, but also from a theoretical
point of view %standpoint
this scenario %state-of-affairs (nmg) 
poses challenging problems. 
Several extensions are being developed; 
%mainly from an analytical point of view; %presented; 
some remarkable ones
are related to the memristive systems of Chua and Kang
\cite{chuakang} and
the mem-devices (memcapacitors and meminductors) recently introduced
by Di Ventra {\em et al.}\ \cite{diventra09}.
%, or the $\sigma$-$\rho$ devices
%(cf.\ subsection ... below) 
%postulated by Biolek {\em et al.}\ \cite{biolek}.
Somehow, the fundamentals of nonlinear circuit theory are affected, since
its concepts and models seem to be moving beyond their
classical limits, and several questions arise. One may wonder about the limits to which this
framework may extend, but also how to accommodate these 
new devices in a comprehensive taxonomy, and which aspects should
articulate such a taxonomy.
%What are the fundamentals... What is the core
In the present paper we will try to answer some of these questions.

%Trying to answer some of these questions, 
In this regard we will find it useful to describe
the memristor (say, in a charge-controlled setting) not by means of its
flux-charge relation 
%\begin{equation} 
%\label{phiq}
$\varphi=\phi(q)$ 
%\end{equation}
but, instead, via the differentiated
relation $\varphi'=\phi'(q)q'$, that is,
\begin{equation} \label{vMqi}
v = M(q) i,
\end{equation}
where $M(q)=\phi'(q)$ is the so-called {\em memristance.}
The reason to do so is that, in the description of the dynamics of memristive
circuits, there is no reason to keep track of both the flux and the charge
as dynamic variables since their values are constrained
%related 
by the relation
%restriction
$\varphi=\phi(q)$. 
%(\ref{phiq}). 
The memristor has order (or, to avoid terminological ambiguities, {\em
  state order}) one. This means that
every memristor introduces one, but not two, degree(s) of dynamic
freedom, associated with the charge $q$;
describing the device characteristic 
via (\ref{vMqi}) we get rid of the flux %linkage
$\varphi$. A detailed discussion can be found in %Section \ref{sec-dev}
subsection \ref{subsec-mem}.

Now, the key remark in our discussion is to look at (\ref{vMqi}) 
just as a relation involving {\em three} out of the four variables
$q,$ $\varphi,$ $i,$ $v$ (all but $\varphi$), 
regardless of its specific form; actually, as in \cite{fully}
we will work with so-called {\em $q$-memristors,} 
defined (in a time-invariant setting) by a general, fully nonlinear characteristic of
the form
$v  = \eta(q, i)$. 

The
corresponding relation for $\varphi$-memristors is
$i=\zeta(\varphi,v)$ and, in particular, 
the characteristic of Chua's flux-controlled memristors
\cite{chuamemristor71} 
%; these are termed voltage-controlled memristors in ... (not clear))
reads as
\begin{equation} \label{iWv}
i = W(\varphi) v,
\end{equation}
in terms of the {\em memductance} $W(\varphi)$; the identity (\ref{iWv}) again
involves three out of the four variables mentioned above (in this case, all
but $q$). 

\ %\vspace{3mm}

\noindent {\bf Memcapacitors and meminductors.} In a natural way, this makes us wonder about 
the remaining relations involving the other two combinations of three
variables, namely, $q,$ $\varphi,$ $v$ and $q,$ $\varphi,$ $i$. This
will lead to the (voltage-controlled) {\em memcapacitor}, with a characteristic of the form
\begin{eqnarray}
q & = & C_{\mathrm m}(\varphi)v, 
\end{eqnarray}
and the (current-controlled) {\em meminductor}, for which,
\begin{eqnarray}
\varphi & = & L_{\mathrm m}(q) i. 
\end{eqnarray}
These devices have been introduced by Di Ventra {\em et al.}\
in \cite{diventra09} and are discussed in
subsection \ref{subsec-memdev}, using again more general, fully nonlinear characteristics.

\ %\vspace{3mm}

\noindent {\bf Hybrid memristors.} Finally, for the sake of completeness
one may think about a characteristic actually involving all {\em four}
variables $q,$ $\varphi,$ $i,$ $v$. For reasons
detailed later, we will consider two different settings, defined by the
fully nonlinear characteristics
%\begin{equation}
$v =  \psi(q, \varphi, i)$ 
%\end{equation}
and
%\begin{equation}
$i =  \xi(q, \varphi, v).$
%\end{equation}
In particular, when these relations are linear in $i$ and $v$,
respectively, we get
\begin{equation}
v =  M_{\mathrm h}(q, \varphi) i 
\end{equation}
and
\begin{equation}
i =  W_{\mathrm h}(q, \varphi) v.
\end{equation}
We will call these circuit elements {\em hybrid memristors}, since their
memristance and memductance involve both the charge $q$ and the flux $\varphi$.
%It is worth emphasizing that t
These novel circuit elements,
together with their fully nonlinear variants, are
introduced here for mathematical completeness, but
they also provide a natural framework to accommodate
actual physical devices in which 
different memory effects (namely, memristive, memcapacitive and/or
meminductive features) coexist.
Find details in subsections \ref{subsec-hybrid} and \ref{subsec-examples}.

%might describe physical devices or provide more accurate models
%for already existing ones.

%Notably, 
%the use of a fully implicit form $H(q, \varphi)=0$ leads to ... instead of...

\ %\vspace{3mm}

These relations comprise all possible 
combinations of the four variables $q$, $\varphi$, $i$, $v$. This
means that the aforementioned devices and their time-varying analogs 
exhaust what  may be called 
{\em first order circuit theory}, within the setting defined by
items (a), (b) and (c) on p.\ \pageref{page-abc}. Here the term
``order'' is used to mean {\em differential order}
(cf.\ %Definition \ref{defin-orderone} in 
Section \ref{sec-dev}) and, accordingly,
the 
expression ``first order circuit'' reflects the fact that only the branch
currents and voltages $i$, $v$ and their first integrals $q$,
$\varphi$ are involved (compare Definitions \ref{defin-orderone} and
\ref{defin-ordertwo} in Sections \ref{sec-dev} and
\ref{sec-other}).
Certainly, the reader should not misunderstand
our use of this expression with that in elementary
circuit theory, referring to a circuit with only one reactive
element.  %Section 2....

%The (a), (b), (c) setting above 
%implicitly assumes that the devices have order one.

%Number of variables involved in the characteristic; number of dynamic variables

Beyond this framework, {\em higher order devices} involve additional
variables such as the first integrals of the charge and the flux
(second integrals of the branch current and voltage),
%these first integrals will be 
to be denoted by $\sigma$ and $\rho$, respectively.
With the use of these additional variables we may accommodate devices
such as charge-controlled memcapacitors and flux-controlled
meminductors \cite{diventra09}.
%, or the $\sigma$-$\rho$ devices discussed in \cite{biolek}. 
Another
extension of the theory is related to the {\em mem-systems} introduced
by Chua and Kang \cite{chuakang}. Roughly speaking, the restriction on
the form of the differential relations considered in (a) is no longer
assumed in these systems; by contrast, new variables without a
physical meaning in classical circuit-theoretic terms are allowed, 
providing dynamical states on which the resistance, capacitance and
inductance
(or better, the memristance, memcapacitance and meminductance)
may depend. Higher order devices and mem-systems are briefly discussed 
in Section \ref{sec-other}: a deeper analysis of them is in the scope of future research.
%Kirchhoff laws are the core. Allow for new variables, with a classical
%circuit-theoretical meaning in higher order devices.
%And remove the restrictions stated in (a) in mem-systems, introducing
%new variables which do not necessarily correspond to currents,
%voltages or their integrals.

The use of these types of devices in applications, as
well as the numerical simulation of the dynamics of nonlinear circuits including
them, requires some effort also at the circuit modelling level. 
We will consider circuit modelling
aspects in Section \ref{sec-nodal}. Specifically, we 
analyze there the differential-algebraic models arising from the nodal analysis of
nonlinear circuits including all possible types of first order
devices.
Our goal in this regard will be the characterization of the
{\em index} of these models, a notion 
which defines many key properties of differential-algebraic
systems
\cite{bre1, grima, hai2, kmbook, rabrheintheo, wsbook}.
%We will briefly discuss second order devices and mem-system in
%Section ***. 
This way we will arrive at a general index characterization (cf.\
Theorems \ref{th-i1} and \ref{th-i2} in Section \ref{sec-nodal})
from which the results discussed in
\cite{et00, wsbook, memristors, carentop, carenhabil} can be derived
as particular cases which apply to circuits including only
certain types of devices with differential order one; another 
particular case of importance concerns circuits with memcapacitors and
meminductors,
whose index is so far unexplored in the literature.
We will also extend the notion of a topologically
degenerate configuration to this broader setting.

%for ``classical'' nonlinear circuits, as well as those in \cite{memristor}
%** Finally, Section \ref{sec-con} compiles some
%concluding remarks. ***

\section{First order devices}
\label{sec-dev}

\begin{defin} \label{defin-orderone}
A circuit device is said to have differential order one
%{\em first order device} is a circuit element 
if it is defined by a
$C^1$-characteristic of the form
\begin{equation}
h(q, \varphi, i, v, t)  =  0,
\end{equation}
where at least one of the partial derivatives %$\partial h / \partial q$
$h_q$, $h_{\varphi}$
does not vanish identically.
%with
%\begin{equation}
%\left( \frac{\partial h}{\partial q}, \frac{\partial h}{\partial
%  \varphi} \right) \not\equiv 0.
%\end{equation}
\end{defin}

Note that both $h_q$ and $h_{\varphi}$ may of course vanish 
for specific values of the circuit variables. Definition \ref{defin-orderone} 
is supported on the fact that when both $h_q$
and $h_{\varphi}$ vanish identically, the only branch variables involved
are the current and/or the voltage. This amounts to what may be called
(differential) order zero devices, not involving any dynamics. These are
%\begin{itemize}
%\item 
nonlinear resistors, %where we use this term to mean 
and voltage and current sources.
%\end{itemize}
%\noindent 
By a nonlinear resistor we mean any device
relating current and voltage in an algebraic (i.e., non-differential)
manner, such as a diode.

The use of the term ``order'' is worth a digression. This term
has different senses in mathematics, and at least two 
apply in our context. On the one hand it means the order of derivation of
a system of differential equations; in this sense, Newton's law yields a second 
order system. We will use the expression {\em differential order}
to mean this. On the other hand, ``order'' is also
%is aimed at making a difference
%with the meaning of the term ``order'' 
used to mean the number of dynamic variables in a system of differential
equations; the term is often used in this sense in circuit theory, %in order 
to refer to the number of state variables 
associated with a given device. We will describe the latter
as the {\em state order} of a device. Throughout the document and when no
label is used, by ``order'' we refer to the differential order; e.g.\
a first order device is a device with differential order one.

\vspace{3mm}

As detailed in the sequel, % forthcoming subsections, 
the list of devices with
differential order one will include
\begin{itemize}
\item %state order one devices:
capacitors, inductors, $q$-memristors and $\varphi$-memristors, all of
them with state order one; %and
\item %state order two devices: 
voltage-controlled memcapacitors, current-controlled meminductors, and
hybrid memristors, with state order two.
\end{itemize}
%Notation derivatives; %notation $\not\equiv$

\noindent 
In practice, at least %$h$ has at least two derivatives which do not
%vanish identically.
%o We implicitly assume that one of the other three... o arriba ``at
%least 
two of the derivatives 
$h_q$, $h_{\varphi}$, $h_i$, $h_v$ 
%and in particular at least one of 
%$h_q$, $h_{\varphi}$, does 
do not vanish identically in first order devices.
It is also worth noting that in Definition \ref{defin-orderone} we implicitly assume 
 the arguments of %the map 
$h$ to be
one-dimensional; however, allowing $q$, $\varphi$, $i$ and $v$ to take
vector values, 
this definition easily accommodates coupled devices.
In this case the non-vanishing requirement on
$h_q$ and/or $h_{\varphi}$ must be replaced by the non-singularity of
the corresponding matrices of partial derivatives. This remark 
will apply throughout the paper, often without explicit mention.

\subsection{%Devices with state order one; 
$q$- and $\varphi$-memristors}
\label{subsec-mem}
%Note that the pairs of variables $q-i$
%and $\varphi-v$ are %already 
%linked by t
The electromagnetic relations $q'=i$ and $\varphi'=v$ link the 
pairs of fundamental circuit variables $q-i$ and $\varphi-v$. First order devices involving the other 
combinations of two out of the four variables $q$,
$\varphi$, $i$, $v$ are the {\em capacitor}, relating $q$ and $v$,
the {\em inductor}, which involves $\varphi$ and $i$, 
and Chua's {\em memristor}, whose characteristic relates $q$ and $\varphi$. %\cite{chuamemristor71}.

The complete dynamical description of the first two devices requires
including the differential relation $q'=i$ for the capacitor (note that
$\varphi$ does not appear in its constitutive relation),
and $\varphi'=v$ for the inductor ($q$ not being involved). 
%Things are different for 
Regarding the memristor, even though both $q$ and $\varphi$ arise
in its characteristic, these variables are linked together and
therefore only one of them introduces a dynamical degree of freedom in
the circuit. For this reason, it will be preferred to get rid of
either the flux or the charge by means of a relation formulated in
terms of the
other three variables, as detailed below. 
This way, not only the capacitor and the inductor but also $q$- and
$\varphi$-memristors will have state order one.

\begin{defin} \label{defin-qmem}
A {\em $q$-memristor} %(resp.\ a $\varphi$-memristor) 
is a device with differential order one, governed by the relations\hspace{-2mm}
\begin{subequations} \label{qmem}
\begin{eqnarray}
q' & = & i \\
v & = & \eta(q, i, t), \label{qmemb}
\end{eqnarray}
\end{subequations}
where $\eta$ is a $C^1$-map for which neither of the derivatives
$\eta_q$, $\eta_i$ vanishes identically.
\end{defin}
When $\eta$ is time-invariant and linear in $i$, we get Chua's
characteristic $v=M(q)i$, where $M(q)$ is the {\em memristance} \cite{chuamemristor71}. The
identity $q(t)=\int^{t}_{-\infty}i(\tau)d\tau$ shows that this
voltage-current relation keeps track of the device history; for this
reason Chua proposed the name {\em memory-resistor}, or {\em memristor}
for short. In the literature, this device is said to be {\em current-controlled} 
but also {\em charge-controlled}, because of the form of the map
$\varphi=\phi(q)$ whose time derivative yields the voltage-current
relation $v=M(q)i$.

For fully nonlinear memristors of the form (\ref{qmem}) the {\em
  incremental memristance} is defined as the derivative $\eta_i(q, i,
t)$ \cite{fully}; the device is called strictly locally passive if
  $\eta_i(q, i, t)>0$ for all $(q, i, t)$ (or if the {\em matrix}
  $\eta_i$ is positive definite in coupled cases).
The requirement that $\eta_i$ does not vanish
identically distinguishes the device from a nonlinear capacitor.
In turn, the non-vanishing condition on $\eta_q$
makes this device actually different from a nonlinear,
current-controlled resistor. Note also that the fully nonlinear form
  (\ref{qmem}) makes it possible to accommodate devices displaying
memristive effects but whose characteristic does not arise as the time
  derivative of a $\varphi$-$q$ relation, contrary to Chua's memristor.

%$\varphi$ is decoupled from the remaining circuit
%equations: state order one; derivative may vanish at certain points;

%Eliminate one variable, but also allow, via a fully nonlinear relation
%\cite{fully}, for the description of devices
%which may display memristive effects not coming from a $\varphi$-$q$ relation.

\begin{defin} \label{defin-phimem}
A {\em $\varphi$-memristor}
is a device with differential order one,  governed by %the relations
\begin{subequations}
\begin{eqnarray}
\varphi' & = & v \\
i & = & \zeta(\varphi, v, t), \label{phimemb}
\end{eqnarray}
\end{subequations}
where $\zeta$ is a $C^1$-map such that neither of the derivatives
$\zeta_{\varphi}$, $\zeta_v$ vanishes identically.
\end{defin}

A time-invariant $\varphi$-memristor for which (\ref{phimemb}) is linear in
$v$ (hence reading as $i=W(\varphi)v$) 
amounts to Chua's flux-controlled memristor \cite{chuamemristor71},
$W(\varphi)$ being the {\em memductance}. In
general, the {\em incremental memductance} is the derivative
$\zeta_v(\varphi, v, t)$, and the device is said to be strictly
locally passive if this derivative is always positive. Again, the non-vanishing requirements on 
$\zeta_{\varphi}$ and  $\zeta_v$  make this device different from a
voltage-controlled resistor and an inductor, respectively.

Nonlinear, $C^1$ circuits composed of (current- and
voltage-controlled) 
resistors, voltage and current sources, capacitors, inductors, and
$q$- and $\varphi$-memristors display the property that
neither the differential nor the state order of their
devices exceed one.
Several dynamical properties of this type of circuits are 
discussed in \cite{fully}. In the next subsections we discuss the main
features of devices with differential order one and state order two.

\subsection{Memcapacitors and meminductors}
\label{subsec-memdev}

The remaining characteristics involving three out of the four
variables $q$, $\varphi$, $i$, $v$ naturally lead to the 
voltage-controlled memcapacitors 
and current-controlled meminductor discussed below. 
%The presence of
%both $q$ and $\varphi$ in the description of these devices means that
Both devices will have state order two.
We consider 
fully nonlinear characteristics for memcapacitors and meminductors; when the maps $\omega$
and $\theta$ below are time-independent and linear in $v$ and $i$,
respectively, we will get 
the circuit elements introduced by Di Ventra {\em et al.}\
in \cite{diventra09}. %, as detailed in the sequel. 

\begin{defin} \label{defin-memcap}
A {\em voltage-controlled memcapacitor}
is a device with differential order one, defined by 
\begin{subequations} \label{memcap}
\begin{eqnarray}
q' & = & i\\
\varphi' & = & v \\
q & = & \omega(\varphi, v, t), \label{memcapc}
\end{eqnarray}
\end{subequations}
where $\omega$ is a $C^1$-map and neither of the derivatives $\omega_{\varphi}$, $\omega_v$ vanishes
identically.
\end{defin}
When $\omega$ is time-invariant and linear in $v$, (\ref{memcapc}) reads as
$q=C_{\mathrm m}(\varphi)v$, where $C_{\mathrm m}$ is the 
{\em memcapacitance} \cite{diventra09}.
The distinct feature of this device is that the memcapacitance depends
on the state variable $\varphi(t) = \int_{-\infty}^t v(\tau)d\tau$,
so that the relation $q(t)=C_{\mathrm m}(\int_{-\infty}^t v(\tau)d\tau)v(t)$
reflects the device history. Be aware of the 
circuit-theoretic meaning of this state variable; when the
memcapacitance is allowed to depend on abstract state variables we are
led to the more general setting of {\em memcapacitive systems} (see \cite{diventra09}
and Section \ref{sec-other} below). It is also worth mentioning that 
the relation $q=C_{\mathrm m}(\varphi)v$ arises as the derivative of
a characteristic $\sigma = \mu(\varphi)$, where $\sigma$ is the
time integral of $q$. This
yields $\sigma'=q=\mu'(\varphi)\varphi'=C_{\mathrm
  m}(\varphi)v$. Notably, the device 
can be described without recourse to 
the (second order) variable $\sigma$ (cf.\ (\ref{memcap})), in contrast to the
charge-controlled memcapacitors discussed in Section \ref{sec-other}.

For an arbitrary $C^1$-map $\omega$, the {\em incremental
memcapacitance} $C_{\mathrm m}$
is defined as the derivative $\omega_v$ and, in general, depends on
$(\varphi, v, t)$. Such a device need not come from a
$\sigma-\varphi$ relation, and 
the requirement that neither $\omega_{\varphi}$ nor $\omega_v$ vanishes
identically makes it actually different from a capacitor or
a memristor, respectively. An instance of a fully nonlinear
memcapacitor arising in a Josephson junction model
can be found in subsection \ref{subsec-examples}.

%\vspace{3mm}

\begin{defin} \label{defin-memind}
A {\em current-controlled meminductor}
is a device with differential order one, governed by
\begin{subequations} \label{memind}
\begin{eqnarray}
q' & = & i\\
\varphi' & = & v \\
\varphi & = & \theta(q, i, t), \label{memindc}
\end{eqnarray}
\end{subequations}
where $\theta$ is a $C^1$-map for which neither of the derivatives $\theta_q$, $\theta_i$ vanishes
identically.
\end{defin}

Again, when the map $\theta$ is linear in $i$ and does not depend on $t$, we get
the characteristic $\varphi = L_{\mathrm m}(q)i = 
L_{\mathrm m}(\int i(\tau)d\tau)i$ considered by 
Di Ventra {\em et al.}\ in \cite{diventra09}. Such a characteristic
can be obtained as the time derivative of a $C^1$-relation $\rho=\kappa(q)$, where $\rho$
is the time integral of $\varphi$; note however that the description
(\ref{memind}) does not involve $\rho$ (compare with the
flux-controlled meminductors in Section \ref{sec-other}).
Now
$L_{\mathrm m}$ is the {\em meminductance}. In general, 
the derivative $\theta_i(q, i, t)$ is the {\em incremental
  meminductance}. The non-vanishing of this derivative makes the
meminductor different from a memristor and, similarly, the assumption that 
the partial derivative
$\theta_q$ does not vanish identically makes the device different from an inductor.

\vspace{3mm}

The appearance of both $q$ and $\varphi$ in the characteristics
(\ref{memcapc}) and (\ref{memindc}) imply that both variables must be
present in the dynamical description of memcapacitors and
meminductors and, therefore, that these devices have state order two; 
this is an important difference with capacitors,
inductors and $q$- and $\varphi$-memristors, for which one dynamic
variable  suffices  (together with $i$ and $v$) to
describe the device behavior.

 \subsection{Hybrid memristors} %: the first-order picture is complete}
\label{subsec-hybrid}

In light of the characteristics (\ref{qmemb}), (\ref{phimemb}),
(\ref{memcapc}) and (\ref{memindc}), from a
mathematical
point of view 
it is somehow natural to complete the picture by considering a
relation which involves all four variables $q$, $\varphi$, $i$ and
$v$. We present below two different settings (dual to each other) in which these four
variables might actually arise. 

These (so-called {\em hybrid}) memristors may account for physical
devices in which memory effects of different nature coexist. The
simultaneous appearance 
of memristive, memcapacitive and/or meminductive phenomena
has been discussed by different authors. For instance, the coexistence
of memristive and memcapacitive effects has been reported to follow
from the formation of local dipoles in nanoscale resistors
\cite{diventra09}, and may also arise in metal-insulator-metal thin films having
thickness between the nanometer and the micrometer scales \cite{mouttet10}.
%In addition to their intrinsically 
%Besides 
%nonlinear inductive nature,
Capacitive and memristive effects coexist with 
the nonlinear inductive nature of a Josephson junction
%are exhibited 
in accurate models of this device \cite{jeltsema10}.
From a modelling point of view, hybrid memristors also
allow for simplified descriptions of combinations
of more basic devices. Find examples in subsection
\ref{subsec-examples} below.

\begin{defin} A {\em current-controlled hybrid memristor}
is a device with differential order one, defined by the relations
\begin{subequations}
\begin{eqnarray}
q' & = & i\\
\varphi' & = & v \\
v & = & \psi(q, \varphi, i, t), \label{hybc}
\end{eqnarray}
\end{subequations}
where $\psi$ is a $C^1$-map and none of the derivatives $\psi_q$, 
$\psi_{\varphi}$, $\psi_i$ vanishes
identically.
\end{defin}

The terminology is motivated by the fact that, in cases in which
(\ref{hybc}) is time-independent and linear in $i$, this relation takes the form 
%\begin{equation}
$v =  M_{\mathrm h}(q, \varphi) i$, 
%\end{equation}
providing an analog of Chua's memristor in which 
the so-called {\em hybrid memristance} $M_{\mathrm h}(q, \varphi)$ now
depends on both the charge and the flux; both variables introduce
memory effects on the device, because of the relations
$q(t)=\int_{-\infty}^{t}i(\tau)d\tau$,
$\varphi(t)=\int_{-\infty}^{t}v(\tau)d\tau$. 
%Not (maybe) $M_\mathrm{h}$. 
In general, the {\em incremental hybrid
memristance} is defined as the derivative $\psi_i$, and may depend
not only on $q$ and $\varphi$ but also on the current $i$ and on $t$.
The hybrid memristor is said to be strictly locally passive when the
incremental hybrid memristance is positive (or positive definite in coupled cases).
%%Model as memristive systems... not so clear

%Memcapacitive and meminductive effects?

\begin{defin} A {\em voltage-controlled hybrid memristor}
is a device with differential order one, governed by 
\begin{subequations}
\begin{eqnarray}
q' & = & i\\
\varphi' & = & v \\
i & = & \xi(q, \varphi, v, t), \label{hybbisc}
\end{eqnarray}
\end{subequations}
where $\xi$ is a $C^1$-map for which none of the derivatives $\xi_q$, 
$\xi_{\varphi}$, $\xi_v$ vanishes
identically.
\end{defin}

%\begin{equation}
Cases in which $\xi$ is linear in $v$ and time-invariant yield
$i =  W_{\mathrm h}(q, \varphi) v$, where $W_{\mathrm h}$ is now the {\em hybrid
  memductance}. In this situation, when $W_{\mathrm h}$ is non-singular for all
values of $q, \varphi$, then obviously the device
can be recast in a current-controlled form and vice-versa.
Again, in general the {\em incremental hybrid
  memductance} is defined by the derivative  $\xi_v(q, \varphi, v,
t)$ and, as before,
the device is said to be strictly locally passive when the
incremental hybrid memductance is positive (or positive definite).

%\end{equation}

\vspace{3mm}

The non-vanishing requirements on the derivatives within the
characteristics (\ref{hybc}), (\ref{hybbisc})
distinguish these devices from the $q$- and $\varphi$-memristors,
memcapacitors and meminductors discussed above. Actually, the
non-vanishing of these derivatives implies that, at least locally,
both the charge and the flux can be written in terms of the remaining
circuit variables. Focusing the attention on the current-controlled
case, the fact that $\psi_q \neq 0$ makes it possible to recast (\ref{hybc}), via
the implicit function theorem, as
\begin{equation} q = \alpha(\varphi, i, v, t), \label{hybq} \end{equation}
and the device exhibits a (generalized) memcapacitance
%\[ \left. \frac{\partial \alpha}{\partial v}\right|_{(\varphi, i, v, t)} = 
%\left. \left( \frac{\partial \psi}{\partial q} \right) \right|^{-1}_{(\alpha(\varphi, i,
%v, t),\varphi, i, t)}\]
\[ \alpha_v (\varphi, i, v, t) = 
\psi_q^{-1}(\alpha(\varphi, i,
v, t),\varphi, i, t).\]
Similarly, the non-zero nature of the partial derivative
$\psi_{\varphi}$ yields
\begin{equation} \varphi = \beta(q, i, v, t), \label{hybphi} \end{equation}
for some locally defined map $\beta$, the (generalized) meminductance being
%\[ \left. \frac{\partial \alpha}{\partial v}\right|_{(\varphi, i, v, t)} = 
%\left. \left( \frac{\partial \psi}{\partial q} \right) \right|^{-1}_{(\alpha(\varphi, i,
%v, t),\varphi, i, t)}\]
\[ \beta_i (q, i, v, t) = 
-\psi_{\varphi}^{-1}(q, \beta(q, i, v, t), i, t) 
\psi_i(q, \beta(q, i, v, t), i, t), 
\]
as a consequence of the implicit function theorem. Similar remarks
apply to voltage-controlled hybrid memristors. 

Noteworthy, removing the non-vanishing requirements on the derivatives,
the relations (\ref{hybc}), (\ref{hybbisc}), (\ref{hybq}) and (\ref{hybphi}) 
account for the characteristics %constitutive relations 
of all previous devices.
From a dynamic point
of view, a remarkable difference between hybrid memristors and both
$q$- and $\varphi$-memristors is that hybrid ones have state order
two, since both the charge and the flux are needed 
%necessarily involved 
in the device description. This implies, for instance,
that the dynamics of a circuit with one hybrid memristor and without reactive
elements lies on the plane.
%, allowing e.g.\ for oscillatory effects
%which cannot be displayed by a system with just one ($q$- or
%$\varphi$-) memristor in the absence of reactive elements. 
The
actual dynamical properties of circuits with hybrid memristors are in
the scope of future research. 
%, which in any case should be driven by
%actual modelling advantages or potential applications of these devices.

%Motivate the following?: no difficulties from an analytical point of
%view, in the sense that in topologically non-degenerate settings a
%state space reduction is feasible... [More]

%\

%$\lambda$, 
%$\kappa$, $\iota$, $\upsilon$, still free]

%In particular, $H(q, \varphi)=0$ yields
%\[H_q(q, \varphi)i + H_{\varphi}(q, \varphi)v=0\] and the non-singularity of
%... or ... leads to ... / ...

%Term also $(q, \varphi)$-memristors

\subsection{Examples}
\label{subsec-examples}

\vspace{3mm}

\begin{figure}[ht]
\centering
\ \mbox{} \hspace{3mm} \
\parbox{2.5in}{
\epsfig{figure=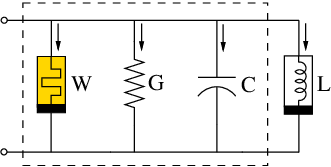, width=0.47\textwidth}
}
\ \mbox{} \hspace{26mm} \
\parbox{2.5in}{\vspace{1mm}
\epsfig{figure=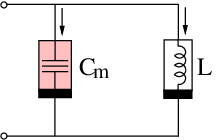, width=0.30\textwidth}
}\hspace{-10mm}
\vspace{2mm}
\caption{(a) Josephson junction. (b) Equivalent circuit
  with a fully nonlinear memcapacitor.}
\label{fig-joseph}
\end{figure}

\noindent {\bf Example 1.} Our first example attempts to illustrate how the
use of fully nonlinear characteristics %may 
allows for simplified
device descriptions. %of certain devices. 
Specifically, we will provide
an accurate description of a Josephson junction by means of a fully
nonlinear memcapacitor, which accounts for all parasitic effects,
connected in parallel to a nonlinear inductor (cf.\ Figure
\ref{fig-joseph}).

\vspace{3mm}

As detailed in \cite{jeltsema10}, realistic models of a Josephson
junction should take into not only the usual nonlinear inductive
relation
$i_l = I_0 \sin (k_0 \varphi_l)$
for certain physical constants $I_0$, $k_0$ (see e.g.\
\cite{chuadesoerkuh}), but also the presence 
of memristive, resistive and capacitive effects; an accurate
equivalent circuit of the Josephson junction is defined by the
parallel connection of these four elements (cf.\ \cite{jeltsema10}),
as depicted in Figure \ref{fig-joseph}(a).

The device on the left of Figure \ref{fig-joseph}(a) is a
%flux-controlled 
$\varphi$-memristor of Chua type, which as reported in \cite{jeltsema10}
captures the presence of
a small current component given by
\begin{equation*}
i_w = I_1 \cos (k_1 \varphi_w) v, %\label{josephw}
\end{equation*}
for certain constants $I_1$, $k_1$; here $v$ is the port voltage.

A linear resistor and a linear capacitor in parallel
are also present in the description provided in \cite{jeltsema10},
being denoted by $G$, $C$, respectively, in Figure
\ref{fig-joseph}(a). 
%We may even accommodate nonlinear effects in
%these two parasitic elements and assume that they are defined by the
%relations
%\begin{equation}i_g = \gamma_g(v)\label{josephg}
%\end{equation}
%and 
%\begin{equation}q_c = \gamma_c(v),\label{josephc}
%\end{equation}
%respectively. In the linear case these identities amount to 
These elements are defined by the relations $i_g = Gv$
and $q_c = Cv$.

Now, the parallel connection of the memristor and the resistor is
obviously governed by the current-voltage relation 
%\begin{equation*}
$i_{wg} = I_1 \cos (k_1 \varphi_w) v + Gv$ %\label{josephwg}
%\end{equation*}
or, equivalently, by a charge-flux characteristic of the form
\begin{equation}
q_{wg} = %\frac{I_1}{k_1} 
(I_1/k_1) \sin (k_1 \varphi_w)  + G\varphi_w, \label{qjosephwg}
\end{equation}
%where $\varphi_{wg}$ 
where we use the fact that the memristor flux %$\varphi_w$
is the time-integral of the port voltage $v$. The expression depicted
in (\ref{qjosephwg}) shows that the
parallel connection of the memristor and the resistor is itself a
%flux-controlled 
$\varphi$-memristor. %of Chua type.

In turn, the parallel connection of the original memristor, the resistor and
the capacitor can be described as a single device by setting $q=q_{wg}
+ q_c$, $\varphi=\varphi_w$. Indeed, denoting by $i$ the sum of the currents
%flowing 
through the memristor, the resistor and the capacitor,  we
get
\begin{subequations}
\begin{eqnarray}
q' & = & i \\
\varphi' & = & v \\
q & = & (I_1/k_1) \sin (k_1 \varphi)  + G\varphi + Cv. \label{josephmemcap}
\end{eqnarray}
\end{subequations}
This corresponds to a %fully nonlinear, 
time-invariant, voltage-controlled
memcapacitor for which the constitutive relation 
$q=\omega(\varphi, v)$ in (\ref{memcapc}) takes the specific form depicted in
(\ref{josephmemcap}). The corresponding equivalent circuit for the
Josephson junction is displayed
in Figure \ref{fig-joseph}(b).

\vspace{2mm}

\noindent {\bf Example 2.} Suppose that the resistor
within the Josephson junction model depicted in Figure \ref{fig-joseph}(a)
exhibits a memristive effect which
makes the conductance $G$ depend on $q$. As detailed in what follows,
the resulting parallel connection of a flux-controlled and a
charge-controlled memristor would yield a simple instance of a 
hybrid memristor.

\begin{figure}[ht]
\centering
\ \mbox{} \hspace{1mm} \
\parbox{2.5in}{
\epsfig{figure=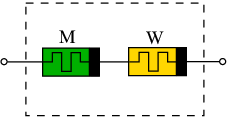, width=0.362\textwidth}
}
\ \mbox{} \hspace{11mm} \
\parbox{2.5in}{\vspace{-0.4mm}
\epsfig{figure=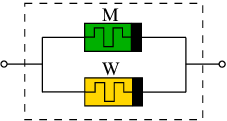, width=0.35\textwidth}
}
\vspace{2mm}
\caption{Examples of (a) a current-controlled and (b) a voltage-controlled hybrid memristor.}
\label{fig-hybrid}
\end{figure}

In broader generality, consider both the series and the 
parallel connection
%Simple instances of hybrid memristors are defined by the connection 
of a charge-controlled and a
flux-controlled memristor of Chua type, as displayed
in Figure \ref{fig-hybrid}. The charge-controlled memristor 
and the
flux-controlled one are painted in green and yellow, respectively.
%In spite of their simplicity, these examples are not trivial: 
As before, the goal is %, again, 
to provide a 
dynamical description of each connection as a two-terminal device in terms of a single set
of variables $q$, $\varphi$, $i$, $v$.

In both cases, the subscripts 1 and 2 will correspond to 
variables associated with the charge-
and the flux-controlled memristor, respectively. The 
charge-controlled memristor is assumed to be governed by a relation of the form
$\varphi_1 = \phi(q_1)$, with memristance $M(q_1)$, and the flux-controlled one
is defined by $q_2 = \gamma(\varphi_2)$, with memductance $W(\varphi_2)$. 
Within Figure \ref{fig-hybrid}(a) %the series 
(resp.\ \ref{fig-hybrid}(b)),
%parallel) configuration,  
the memductance $W(\varphi_2)$ (resp.\ the memristance $M(q_1)$) is assumed not to vanish.

%\vspace{3mm}

%\noindent 
In the series connection of Figure \ref{fig-hybrid}(a),
elementary circuit theory yields $i=i_1 = i_2$, $v=v_1 + v_2$. In order to
arrive at a dynamical description in terms of $i$, $v$ and single
variables $q$, $\varphi$, we set $q=q_1$ and $\varphi=\varphi_1 + \varphi_2$.
Obviously, this yields the relations $q'=i$ and $\varphi'=v$ but, more
important, allows for the description of the voltage-current relation in the form
\begin{eqnarray} v & = & M(q_1)i + (W(\varphi_2))^{-1}i = [M(q) +
  (W(\varphi-\varphi_1))^{-1}]i = [M(q) +
  (W(\varphi-\phi(q_1)))^{-1}]i = \nonumber \\
& = & [M(q) +   (W(\varphi-\phi(q)))^{-1}]i. \nonumber
 \end{eqnarray}
This corresponds to a current-controlled hybrid memristor for which the
characteristic is linear in $i$, the hybrid memristance being
\[M_{\mathrm h}(q, \varphi) = M(q) + (W(\varphi-\phi(q)))^{-1}.
\]
The physical meaning of the variable $q$ is worth some additional remarks.
Because of the identities $q_1' = i_1 = i_2 = q_2'$, the charges 
$q_1$ and $q_2$ differ in a constant which would be fixed by the
initial conditions in Chua's memristors. In turn, $q$ is defined
up to a constant, and therefore can be understood
to describe any of both charges except for a fixed
quantity. Mathematically, 
setting $q=q_1 + k$ for any real constant $k$, we
would get a dynamical description of the device
which,  except 
for an affine change of coordinates, amounts to the one above; 
the description above assumes $k=0$ so that $q$ actually equals the charge
$q_1$.

It is also worth noting that if $M(q) + (W(\varphi-\phi(q)))^{-1}$ does not vanish for all
values of $q$, $\varphi$, the device also admits a voltage-controlled
description. This would be the case, in particular, if the original
Chua %charge-controlled and flux-controlled 
memristors are strictly locally passive
(i.e.\ if $M>0$, $W >0$ everywhere); in this case the hybrid memristor would
itself be strictly locally passive since $M_{\mathrm h}$ would be strictly positive.

%\vspace{1.5mm}

Setting $q=q_1 + q_2,$ $\varphi=\varphi_2$, the reader can proceed analogously
in order to describe the parallel configuration of Figure \ref{fig-hybrid}(b)
as a voltage-controlled hybrid memristor, with a characteristic linear in $v$ and
hybrid memductance
\[W_{\mathrm h}(q, \varphi) = (M(q-\gamma(\varphi)))^{-1} + W(\varphi).\]

%\vspace{1.5mm}

These examples show that, even in simple cases, these devices
pose interesting problems from the modelling point of view.
Needless to say, the scope of hybrid memristors goes however beyond this
type of examples. 
The definition of this circuit element is aimed at modelling 
devices which are not reducible to a simple 
connection of $q$- and $\varphi$-memristors, and which may capture
the coexistence of different memory effects.
In such devices the
characteristics need not be linear in $i$ and $v$, 
and the charge and the flux within the (incremental) memristance or 
memductance might interact in more intricate ways.

%\newpage

\section{Higher order devices and mem-systems}
\label{sec-other}

The devices discussed in Section \ref{sec-dev} virtually fill the scope of
circuit theory within the limits defined by the use of 
the variables $q$, $\varphi$, $i$, $v$ and 
under the restrictions stated in
items (a), (b) and (c) on p.\ \pageref{page-abc}. 
The frontier of nonlinear circuit theory in this setting
is defined by hybrid memristors, which involve all four fundamental
circuit variables. However, other
devices recently introduced %discussed %in the literature
are located beyond these limits: on the one hand, certain circuit elements
involve not only those four variables but also the
time-integrals $\sigma$, $\rho$ of the charge and the flux,
respectively; these include charge-controlled memcapacitors and
flux-controlled meminductors \cite{diventra09}, but also 
%eventual 
devices directly relating $\sigma$ and $\rho$, as proposed
in \cite{biolek}; all of them
would have differential order two.
On the other hand, allowing for %the introduction of 
new state variables, without the restriction imposed by (a), leads to the
so-called {\em mem-systems}, which in general cannot be
accommodated in the framework of Section \ref{sec-dev}. In this Section we briefly 
discuss some modelling aspect of these devices, leaving a deeper analysis
%of these circuit elements 
%defines a line 
for future research.

\subsection{Second order devices}

\begin{defin} \label{defin-ordertwo}
A circuit device is said to have differential order two
%{\em first order device} is a circuit element 
if it is defined by a
$C^1$-characteristic of the form
\begin{equation}
h(\sigma, \rho, q, \varphi, i, v, t)  =  0,
\end{equation}
where at least one of the partial derivatives %$\partial h / \partial q$
$h_{\sigma}$, $h_{\rho}$
does not vanish identically.
%with
%\begin{equation}
%\left( \frac{\partial h}{\partial q}, \frac{\partial h}{\partial
%  \varphi} \right) \not\equiv 0.
%\end{equation}
\end{defin}

The non-vanishing of $h_{\sigma}$ and/or $h_{\rho}$ implies that at
least one of the differential relations $\sigma'=q$ (and in turn $q'=i$)
or $\rho'=\varphi$ (together with $\varphi'=v$)
must be included in the dynamical description of the device. The
differential order two nature of them stems from the identities
$\sigma''=i$ and $\rho''=v$. Akin to differential order one devices,
in practice at least {\em two} of the derivatives 
$h_{\sigma}$, $h_{\rho}$, $h_q$, $h_{\varphi}$, $h_i$, 
$h_v$ will not vanish identically.

\begin{defin} A {\em charge-controlled memcapacitor} is a device with
  differential order two, defined by the relations
\begin{subequations} \label{charmemcap}
\begin{eqnarray}
\sigma' & = & q \\
q' & = & i\\
%\varphi' & = & v \\
 v & = & \nu(\sigma, q, t), \label{charmemcapc}
\end{eqnarray}
\end{subequations}
where $\nu$ is a $C^1$-map for which neither of the derivatives
$\nu_{\sigma}$, $\nu_q$ vanishes identically.
\end{defin}

\begin{defin} A {\em flux-controlled meminductor} is a device with
  differential order two, governed by
\begin{subequations} \label{fluxmemind}
\begin{eqnarray}
\rho' & = & \varphi \\
%q' & = & i\\
\varphi' & = & v \\
i & = & \chi(\rho, \varphi, t)  \label{fluxmemindc}
\end{eqnarray}
\end{subequations}
where $\chi$ is a $C^1$-map such that neither of the derivatives
$\chi_{\rho}$, $\chi_{\varphi}$ vanishes identically.
\end{defin}
From the descriptions (\ref{charmemcap}) and (\ref{fluxmemind}) it
follows that
charge-controlled memcapacitors and flux-controlled
meminductors have state order two; note that neither $\rho$ nor
$\varphi$ (resp.\ $\sigma$, $q$) arise in the
characteristic (\ref{charmemcapc}) (resp.\ (\ref{fluxmemindc})).

In particular, when $\nu$ and $\chi$ are time-invariant and linear in $q$ and
$\varphi$, respectively, one gets the devices introduced by Di Ventra
{\em et al.}\ in \cite{diventra09}, for which the characteristics
(\ref{charmemcapc}) and (\ref{fluxmemindc}) read as
\begin{equation}
v = E_{\mathrm m}(\sigma)q \label{vEmq}
\end{equation}
and 
\begin{equation} 
i = {\mathcal R}_{\mathrm m}(\rho) \varphi. \label{iRmrho}
\end{equation}
Here $E_{\mathrm m}$ and ${\mathcal R}_{\mathrm m}$ are the 
inverse memcapacitance
and the inverse meminductance, which introduce memory effects in
the circuit because of their dependence on $\sigma$ and $\rho$,
respectively. 

Noteworthy, the relations (\ref{vEmq}) and (\ref{iRmrho}) arise as
the differentiated form of certain mappings
$\varphi=\gamma(\sigma)$ and $q=\delta(\rho)$, via the relations
$\sigma'=q$, $\rho'=\varphi$. In a natural way this leads to
other second order devices, such as those relating $\sigma$ and
$\rho$ (cf.\ \cite{biolek}). The analysis of the dynamics of nonlinear
circuits including these or other higher order devices is an open problem,
beyond the scope of the present paper.

\subsection{Mem-systems}

Mem-systems, originally introduced by Chua and Kang in
\cite{chuakang}, are characterized by the removal of the classical
electrical meaning of the state variable which introduces memory into the
different devices. This means that new variables arise in the
model, and that the specific form of the differential relations
within item (a) on p.\ \pageref{page-abc} is no longer assumed.
For the sake of brevity, we restrict the discussion to the mem-systems which generalize the
devices considered in Definitions
\ref{defin-qmem}, \ref{defin-phimem},
\ref{defin-memcap} and \ref{defin-memind}, using again fully nonlinear characteristics.

{\em Memristive systems} are defined either by a system of the form
\begin{subequations}
\begin{eqnarray}
x' & = & f(x,i) \\
v & = & \tilde{\eta}(x, i, t), \label{memsys1b}
\end{eqnarray}
\end{subequations}
or by
\begin{subequations}
\begin{eqnarray}
y' & = & g(y,v) \\
i & = & \tilde{\zeta}(y, v, t). \label{memsys2b}
\end{eqnarray}
\end{subequations}
The distinct feature with respect to $q$- and
$\varphi$-memristors is that now the memristance
and the memductance depend on {\em arbitrary} state variables $x$,
$y$, respectively, which do not even need to be 
one-dimensional; making $x \equiv q$, $f(q,i)=i$ and $y \equiv
\varphi$, $g(\varphi, v)=v$ we get the $q$- and
$\varphi$-memristors introduced in Definitions \ref{defin-qmem} and \ref{defin-phimem}.
Systems in which the characteristics (\ref{memsys1b}) and
(\ref{memsys2b}) 
amount to $v=M(x,i)i$ and $i=W(y,v)v$, respectively, describe
the settings originally discussed by Chua and Kang in
\cite{chuakang}. 

Analogously, a {\em memcapacitive system} is
governed by the system
\begin{subequations}
\begin{eqnarray}
q' & = & i\\
z' & = & h(z,v) \\
q & = & \tilde{\omega}(z, v, t), 
\end{eqnarray}
\end{subequations}
whereas the equations defining a {\em meminductive system} are
\begin{subequations}
\begin{eqnarray}
u' & = & p(u,i) \\
\varphi' & = & v \\
\varphi & = & \tilde{\theta}(u, i, t).
\end{eqnarray}
\end{subequations}
Again, with $\tilde{\omega}(z, v)=C_{\mathrm m}(z)v$, $\tilde{\theta}(u, i)=L_{\mathrm m}(u)i$ we get
the systems analyzed by Di Ventra {\em et al.}\ in
\cite{diventra09}. 

Mem-systems have many potential applications,
not only in electronics; see e.g.\ \cite{pershin09, %pershin10a, pershin10b,
  pershin10c} and references therein. These systems are likely to define an
active field of research in the near future.

%Can be modelled as a memristive system. Not so clear, i, v.

%In all cases, the electrical meaning of the devices 
%in this generality the electrical meaning of the model now relies only on Kirchhoff laws...

\section{Nodal analysis of first order circuits}
\label{sec-nodal}

From a computational point of view, models of the form (\ref{basic})
offer some difficulties for numerical simulation. 
This is due to the fact that an automatic computation of the loop and
cutset matrices $B$, $D$ is difficult to perform in practice, specially
in high scale integration circuits. For this reason, it is often preferred 
to introduce the {\em node potentials} $e$ in the model, and describe
the circuit equations using nodal analysis. This is the case in most
circuit simulation problems, notably in SPICE and its commercial
variants, which set up the circuit equations using Modified Nodal
Analysis (MNA); cf.\ \cite{et00, gun1, gun2, reis10,
wsbook, takamatsu10, carentop, carenhabil}.

As detailed below, the models arising from nodal analysis naturally
take the form of a differential-algebraic equation (DAE) \cite{bre1, grima,
  kmbook, rabrheintheo, wsbook}. The main problem in the analysis of
such differential-algebraic models is the characterization of the {\em
  index}, a concept which, roughly speaking, extends the notion
of the Kronecker-Weierstrass index of a matrix pencil \cite{gantma}
and measures the numerical
difficulties faced in simulation \cite{bre1, hai2}. 
Index one and index two systems
require specific numerical techniques, and because of this it is
important to characterize the circuit configurations which lead to
models with these indices. Therefore, we undertake in this Section the index analysis of
nodal models of general first order circuits, using the tractability
index framework \cite{grima, maeAN, maeanew, carenhabil} and extending 
the results discussed in \cite{et00, wsbook, memristors, carentop} to
circuits with (voltage-controlled)
memcapacitors, (current-controlled) meminductors, and hybrid memristors. Nonlinear
circuits with Chua's memristors, memcapacitors and meminductors define
a particular case of great interest in current applications 
\cite{chuamemristor08, itohchua09, kavehei10, 
%messias10, muthus10, muthuskokate,  
pershin, pershin09, pershin10a, pershin10b,
pershin10c}.

Also from an analytical point of view it is important to characterize
index one circuit configurations. In index one systems, {\em all} the dynamic
variables of the different devices contribute to the state dimension
of the problem; more precisely, in index one cases 
the state dimension (also called the {\em order of complexity}) of a first order circuit equals
the sum of the state orders of the devices with differential order one.
By contrast, in higher index problems which arise from
so-called {\em topologically degenerate} configurations the feasible 
values for these dynamic variables are restricted by algebraic
(non-differential) constraints. In a classical setting the
topologically degenerate configurations are VC-loops (loops just
defined by capacitors and --possibly-- voltage sources) and IL-cutsets
(cutsets just including inductors and 
--possibly-- current sources);
note that V-loops and I-cutset are excluded in well-posed problems. As
a byproduct of our analysis we introduce the notion of a topologically
degenerate configuration for circuits including 
memcapacitors and meminductors.

\subsection{The nodal model}

%Using %this 
%the below-introduced reduced incidence matrix, and a
Assuming a time-invariant setting for
simplicity, 
the nodal model of a general first order circuit reads as
\begin{subequations} \label{nodal}
\begin{eqnarray}
q_c' & = & i_c \label{nodala}\\
q_{mc}' & = & i_{mc} \label{nodalb}\\
\varphi_l' & = & A_l^T e \label{nodalc}\\
\varphi_{ml}' & = & A_{ml}^T e \label{nodald}\\
q_m' & = & i_m \\
q_{ml}' & = & i_{ml} \\
q_{hm}' & = & i_{hm} \\
q_{hw}' & = & i_{hw} \\
\varphi_w' & = & A_w^T e \\
\varphi_{mc}' & = & A_{mc}^T e \\
\varphi_{hm}' & = & A_{hm}^T e \\
\varphi_{hw}' & = & A_{hw}^T e \label{nodall} \\ %\pagebreak
0 & = & A_c i_c +A_{mc}i_{mc}+ A_u i_u + A_l i_l + A_{ml} i_{ml} +
A_g \gamma_g(A_g^T e) + A_w \zeta(\varphi_w, A_w^Te)+   \nonumber\\
&&    + A_r i_r +A_m i_m +A_{hm}i_{hm}+A_{hw}i_{hw} + A_j i_s(t) \ \ \
\ \ \ \ \ \ \ \label{nodalm} \\
0 & = & q_c - \gamma_c(A_c^Te) \label{nodaln}\\
0 & = & q_{mc} - \omega(\varphi_{mc}, A_{mc}^Te) \label{nodalo}\\
0 & = & v_s(t)-A_u^T e \label{nodalp}\\
0 & = & \varphi_l-\gamma_l(i_l) \label{nodalq}\\
0 & = & \varphi_{ml}-\theta(q_{ml}, i_{ml}) \label{nodalr}\\
0 & = & \gamma_r(i_r) - A_r^Te \label{nodals}\\
0 & = & \eta(q_m, i_m) - A_m^T e \\
0 & = & \psi(q_{hm}, \varphi_{hm}, i_{hm}) - A_{hm}^Te \\
0 & = & i_{hw} - \xi(q_{hw}, \varphi_{hw}, A_{hw}^T e). \label{nodalv} 
\end{eqnarray}
\end{subequations}
The two distinct features of nodal models are the use of node
potentials $e$ and the description of 
Kirchhoff laws as $Ai=0$, $v=A^T e$, in terms of the 
reduced incidence matrix $A$. Provided that the circuit is connected and that
a reference node has been chosen, this matrix is defined as $A=(a_{ij})$ with
\begin{eqnarray*}
a_{ij} = \left\{
\begin{array}{rl}
1 & \text{ if branch } j \text{ leaves  node } i \\
-1 & \text{ if branch } j \text{ enters node } i \\
0 & \text{ if branch } j \text{ is not incident with node } i,
\end{array}
\right.
\end{eqnarray*}
for all nodes except for the reference one. In (\ref{nodal}) 
this matrix is %will be
partitioned by columns according to the electrical nature of the
corresponding branches.

Several additional 
remarks are in order. The size of the model (\ref{nodal})
owes to its very general
nature, which accommodates a
great variety of devices, namely, capacitors, inductors,
current-controlled and voltage-controlled resistors, $q$-memristors,
$\varphi$-memristors, memcapacitors, meminductors, 
current-controlled and voltage-controlled hybrid memristors, and
voltage and current sources. The subscripts corresponding to these
devices are $c$, $l$, $r$, $g$, $m$, $w$, $mc$, $ml$, $hm$, $hw$, $u$ and
$j$, respectively. For the sake of simplicity capacitors and inductors
are assumed to be voltage- and current-controlled, respectively 
(cf.\ the maps $\gamma_c$ and $\gamma_l$ in (\ref{nodaln}) and (\ref{nodalq})), 
and
we eliminate the branch currents  $i_g$, $i_w$ of voltage-controlled
resistors and $\varphi$-memristors by means of the maps $\gamma_g$ and
$\zeta$, respectively. The map $\gamma_r$ in (\ref{nodals}) defines the
characteristic of current-controlled resistors, whereas $i_s(t)$ and
$v_s(t)$ are the excitations in the current and voltage sources,
respectively. The maps $\eta$, $\omega$, $\theta$,  $\psi$, $\xi$ are
those arising in the characteristics (\ref{qmemb}), (\ref{memcapc}), 
(\ref{memindc}), (\ref{hybc}) and (\ref{hybbisc}), except for the fact
that they account for the whole sets of  $q$-memristors,
memcapacitors, meminductors and hybrid memristors,
and hence need not be scalar; the same holds, of course, for
$\varphi$-memristors and the map $\zeta$.

For later use, denote by $R$, $G$, $C$, $L$, $M$, $W$, 
$C_{\mathrm m}$, $L_{\mathrm m}$, $M_{\mathrm h}$ and
$W_{\mathrm h}$ the incremental 
resistance, conductance, capacitance, inductance, memristance, memductance, 
memcapacitance, meminductance, hybrid memristance and hybrid memductance
matrices, defined by the derivatives $\gamma_r'$, $\gamma_g'$,
$\gamma_c'$, $\gamma_l'$, $\eta_{i_m}$, $\zeta_{v_w}$, $\omega_{v_{mc}}$,
$\theta_{i_{ml}}$, $\psi_{i_{hm}}$, $\xi_{v_{hw}}$, respectively. 
These matrices need not be diagonal, meaning that full coupling is allowed within each 
of these sets of devices. 
%An important requirement in the analysis is that
%some of these devices are strictly locally passive, so that the 
%corresponding matrices are positive definite: an $n \times n$ matrix $P$ is said
%to be positive definite if $u^T P u > 0$ for all $u \in \mathbb{R}^n-\{0\}$;
%we do not require $P$ to be symmetric.
%In the absence of hybrid memristors...

%Simplified model tipo MNA briefly discussed below... 

\subsection{Topologically nondegenerate configurations and index one models}

%Def topo non-deg;

Semiexplicit differential-algebraic equations
are defined by a system of the form
\begin{subequations} \label{seDAE}
\begin{eqnarray}
x' & = & f(x, y, t) \label{seDAEa}  \\
0 & = & g(x, y, t), \label{seDAEb}
\end{eqnarray}
\end{subequations}
where %$f: \mathbb{R}^{r+p} \to \mathbb{R}^{r}, \ g:
%\mathbb{R}^{r+p} \to \mathbb{R}^{p}$. REGULARITY.
$x \in \mathbb{R}^r$ denotes the {\em differential} or
{\em dynamic} variables,
$y \in \mathbb{R}^p$ stands for the {\em algebraic}
ones,
% REGULARITY.
$f \in C^1(\mathbb{R}^{r+p+1},\ \mathbb{R}^r),$ and
$g \in C^1(\mathbb{R}^{r+p+1},\ \mathbb{R}^p)$. 
The DAE (\ref{seDAE}) is said to be index one around a given $(x^*,
y^*, t^*)$ satisfying (\ref{seDAEb}) if
the matrix of partial derivatives $g_y(x^*, y^*, t^*)$ is invertible.
Often, this non-singularity requirement holds everywhere.
This definition applies in the context of the geometric, differentiation and tractability
indices, and for both analytical and numerical purposes it is
important to characterize index one systems in practice;
detailed discussions about the different index notions
can be found in \cite{bre1, grima, hai2, kmbook,
  rabrheintheo, wsbook}. 

The nodal model 
(\ref{nodal}) has a semiexplicit form, the set of dynamic variables being
\begin{equation} \label{dynvar}
q_c, \ q_{mc}, \ \varphi_l, \ \varphi_{ml}, \
q_m, \ q_{ml}, \ q_{hm}, \ q_{hw}, \ \varphi_w, \ \varphi_{mc}, \
\varphi_{hm}, \ \varphi_{hw},
\end{equation}
whereas the algebraic ones are
\begin{equation} \label{algvar}
e,\  i_c, \ i_{mc}, \ i_u, \ i_l, \ i_{ml}, \ i_r, \ i_m, \ i_{hm},
\ i_{hw}. \end{equation}
We address below
the characterization of index one configurations for 
(\ref{nodal}), 
under the assumption that certain circuit matrices are positive definite; recall that
a given matrix $K$ is positive definite if $u^T K u > 0$ for any
non-vanishing real vector $u$, and that this notion expresses
mathematically a strict passivity requirement on the corresponding devices.
The proof proceeds by showing how the non-singularity of the 
matrix defining index one
configurations 
%(as well as
%the one arising in the index two analysis of subsection ...) 
can be 
reduced to a form already analyzed in the 
context of nonlinear circuits without memristive devices
in \cite{wsbook}. The result stated in Theorem
\ref{th-i1} actually motivates the following definition.

\begin{defin}
A first order circuit is said to be {\em topologically nondegenerate}
if it does not display either loops defined by voltage sources,
capacitors and/or memcapacitors, or cutsets composed of current
sources, inductors and/or meminductors.
\end{defin}
This extends a well-known notion for RLC circuits, for which
topologically nondegenerate configurations preclude
loops defined by voltage sources and/or
capacitors, and cutsets composed of current
sources and/or inductors. Stemming from the work of Bashkow \cite{bashkow57}
in the classical RLC setting,
this provides a way to formulate a state space model of the circuit dynamics
by means on the notion of a {\em proper tree}. Note that in our present
setting neither $q$- or $\varphi$-memristors, nor hybrid ones, introduce
topological degeneracies.

\begin{teorema} \label{th-i1}
Assume that the capacitance 
$C$, 
the memcapacitance $C_{\mathrm m}$, the inductance $L$ and the meminductance 
$L_{\mathrm m}$ are non-singular matrices, and that
the resistance $R$, the conductance $G$, the memristances $M$,
$M_{\mathrm h}$ and the memductances
$W$, $W_{\mathrm h}$  are positive definite.

Then the %nodal 
model (\ref{nodal}) is index one if and only if the circuit is 
topologically nondegenerate.
\end{teorema}

\noindent {\bf Proof.} The matrix of partial derivatives of the right-hand side of
(\ref{nodal}) w.r.t.\ all variables but time, to be denoted by $F$, has the form
\begin{eqnarray} \label{Jac}
F = 
\left(\begin{array}{cc} 
0 & F_{12} \\
F_{21} & F_{22}  
\end{array}\right),
\end{eqnarray}
where the block $F_{22}$ stands for the partial derivatives of the
restrictions
(\ref{nodalm})--(\ref{nodalv}) with respect to the algebraic variables
(\ref{algvar}). The non-singularity of $F_{22}$ characterizes index
one configurations, and this matrix reads as
\begin{eqnarray*}
F_{22} = 
\left(\begin{array}{cccccccccc} 
A_g G A_g^T + A_w W A_w^T & A_c & A_{mc} & A_u & A_l & A_{ml} & A_r & A_m & A_{hm} & A_{hw} \\
-C A_c^T & 0 & 0 & 0 & 0 & 0 & 0 & 0 & 0 & 0 \\
-C_{\mathrm m} A_{mc}^T & 0 & 0 & 0 & 0 & 0 & 0 & 0 & 0 & 0 \\
-A_u^T & 0 & 0 & 0 & 0 & 0 & 0 & 0 & 0 & 0 \\
0 & 0 & 0 & 0 & -L & 0 & 0 & 0 & 0 & 0 \\
0 & 0 & 0 & 0 & 0 & -L_{\mathrm m} & 0 & 0 & 0 & 0  \\
-A_r^T & 0 & 0 & 0 & 0 & 0 & R & 0 & 0 & 0 \\
-A_m^T & 0 & 0 & 0 & 0 & 0 & 0 & M & 0 & 0 \\
-A_{hm}^T & 0 & 0 & 0 & 0 & 0 & 0 & 0 & M_{\mathrm h} & 0 \\
-W_{\mathrm h}A_{hw}^T & 0 & 0 & 0 & 0 & 0 & 0 & 0 & 0 & I_{hw} \\
\end{array}\right),
\end{eqnarray*}
the block $I_{hw}$ being an identity matrix whose size is defined
by the number of voltage-controlled hybrid memristors. The
non-singularity of $C$, $C_{\mathrm m}$, $L$, $L_{\mathrm m}$ obviously reduces the
problem to the characterization of the non-singularity of 
\begin{eqnarray*}
\left(\begin{array}{cccccccc} 
A_g G A_g^T + A_w W A_w^T & A_c & A_{mc} & A_u & A_r & A_m & A_{hm} & A_{hw} \\
-A_c^T & 0 & 0 & 0 & 0 & 0 & 0 & 0 \\
-A_{mc}^T & 0 & 0 & 0 & 0 & 0 & 0 & 0 \\
-A_u^T & 0 & 0 & 0 & 0 & 0 & 0 & 0 \\
-A_r^T & 0 & 0 & 0 & R & 0 & 0 & 0 \\
-A_m^T & 0 & 0 & 0 & 0 & M & 0 & 0 \\
-A_{hm}^T & 0 & 0 & 0 & 0 & 0 & M_{\mathrm h} & 0 \\
-W_{\mathrm h}A_{hw}^T & 0 & 0 & 0 & 0 & 0 & 0 & I_{hw} \\
\end{array}\right)
\end{eqnarray*}
and, by means of a Schur reduction \cite{hor0, wsbook},
 the non-singularity of this matrix amounts to that of
\begin{eqnarray*} 
\left(\begin{array}{cccccccc} 
\hspace{-1mm}A_g G A_g^T \hspace{-0.5mm} + \hspace{-0.5mm} A_w W A_w^T
\hspace{-0.5mm} + \hspace{-0.5mm} 
A_r R^{-1} A_r^T \hspace{-0.5mm} + \hspace{-0.5mm} A_m M^{-1} A_m^T 
\hspace{-0.5mm} + \hspace{-0.5mm} A_{hm} M_{\mathrm h}^{-1}A_{hm}^T 
\hspace{-0.5mm} + \hspace{-0.5mm} A_{hw}W_{\mathrm h} A_{hw}^T
& \hspace{-0.5mm} A_c & \hspace{-0.5mm} A_{mc} & \hspace{-0.5mm} A_u \\
-A_c^T & 0 & 0 & 0 \\
-A_{mc}^T & 0 & 0 & 0 \\
-A_u^T & 0 & 0 & 0
\end{array}\right)\hspace{-0.5mm}.
\end{eqnarray*}
This matrix has the structure arising in \cite[Theorem
  5.1(1)]{wsbook} (cf.\ eq.\ (5.43) there); from this result it follows that, in the present
setting, the non-singularity of this matrix relies on the
absence of loops composed of voltage sources,
capacitors and/or memcapacitors and cutsets defined by current
sources, inductors and/or meminductors, as we aimed to show.

\hfill $\Box$

\subsection{Topological degeneracies: Index two}

In presence of the topologically degenerate configurations discussed
above (i.e.\ loops defined by voltage sources,
capacitors and/or memcapacitors, or cutsets composed of current
sources, inductors and/or meminductors), 
the index one condition for the nodal system (\ref{nodal})
fails. In this situation, for both analytical and numerical purposes 
it is important to characterize whether the
model is index two or not.
The index two notion for a DAE is more intricate than the index one
concept introduced above. Again, the reader is referred to 
\cite{bre1, grima, hai2, kmbook,
  rabrheintheo, wsbook} for different approaches to the index
notion. 

In particular, the {\em tractability index} notion, together with 
the projector-based framework supported on it \cite{grima, maeAN, maeanew,
  mae1, wsbook, lframe, carenhabil},
has been proved to
be a valuable tool in circuit simulation \cite{et00, gun1, gun2,
maeMNA, wsbook, carentop, carenhabil}. In order to introduce this notion, we
look at (\ref{nodal}) as a semilinear problem of the form
\begin{equation} E z' = f(z, t), \label{qlDAE} \end{equation}
where $E$ is a block-diagonal matrix block-diag$\{I,0\}$,
and $z$ joins together the variables denoted by $x$ and $y$ in (\ref{seDAE}).
Consider the {\em matrix pencil} $\lambda E- F$ \cite{gantma}, $F$ being the matrix of
partial derivatives $f_z$. As detailed in the references above, 
the pencil is said to have tractability index one if 
%\begin{equation} 
$E_1 = E - FQ$ 
%\end{equation}
is a non-singular matrix. In turn, if $E_1$ is singular, 
we let $Q_1$ be any projector onto $\ke E_1$, 
and the pencil is said to have tractability index two if 
%\begin{equation}
$E_2 = E_1 - F_1Q_1$
%\label{A2Q}
%\end{equation}
is non-singular, where $F_1=F(I-Q)$. 

Iteratively, this approach provides a
general index notion which can be shown to equal the
Kronecker-Weierstrass (or nilpotency) index of the pencil, being well-suited for
computational purposes. Moreover, this concept can be extended to
nonlinear and/or time-varying settings under suitable
assumptions on the system operators; restricting the attention to
DAEs of the form (\ref{qlDAE}) in an index two context, these assumptions
amount to requiring that $Q_1$ be a {\em continuous} projector onto the
kernel of $E_1(z)$. Supported on these
ideas, we show below that the nodal model (\ref{nodal}) is indeed
index two in the presence of degenerate configurations, under
additional passivity requirements.
Note that in the index two
context the {\em normal tree} method of Bryant \cite{bryant59, bryant62}
applies in order to derive a state space equation.

\begin{teorema} \label{th-i2}
Assume that the capacitance  $C$, 
the memcapacitance $C_{\mathrm m}$, the inductance $L$, the meminductance 
$L_{\mathrm m}$, the resistance $R$, the conductance $G$, the memristances $M$,
$M_{\mathrm h}$ and the memductances
$W$, $W_{\mathrm h}$ are positive definite.

Then the nodal model (\ref{nodal}) has tractability index two in the
presence of topologically degenerate configurations.
\end{teorema}

\noindent {\bf Proof.} As indicated above, the model (\ref{nodal}) has
the form depicted in (\ref{qlDAE}) with $E=$block-diag$\{I,0\}$ and
$F$ being the matrix or partial derivatives displayed in
(\ref{Jac}). Letting $Q$ be a projector
onto the kernel of $E$ with the structure block-diag$\{0, I\}$, we
arrive at
\begin{eqnarray*}
E_1 & = &
\left(\begin{array}{cc} 
I & -F_{12} \\
0 & -F_{22}  
\end{array}\right).
\end{eqnarray*}
Note, incidentally, that this makes it clear why the index one
condition arising in Theorem \ref{th-i1} relies on the non-singularity
of $F_{22}$. As indicated in the proof of Theorem \ref{th-i1}, the
presence of loops defined by voltage sources,
capacitors and/or memcapacitors, and/or cutsets composed of current
sources, inductors and/or meminductors makes $F_{22}$ (and hence
$E_1$) a singular matrix.

In order to construct a projector $Q_1$ onto $\ke E_1$, we need to
detail the structure of $F_{12}$ and $F_{21}$ in (\ref{Jac}). These
blocks are defined by the derivatives of
(\ref{nodala})--(\ref{nodall}) w.r.t.\ the algebraic variables
(\ref{algvar}) and the derivatives of the
restrictions
(\ref{nodalm})--(\ref{nodalv}) w.r.t.\ the dynamic variables
(\ref{dynvar}), respectively, and they read as
\begin{eqnarray*}
F_{12} = 
\left(\begin{array}{cccccccccc} 
0 & I_c & 0 & 0 & 0 & 0 & 0 & 0 & 0 & 0 \\
0 & 0 & I_{mc} & 0 & 0 & 0 & 0 & 0 & 0 & 0 \\
A_l^T & 0 & 0 & 0 & 0 & 0 & 0 & 0 & 0 & 0 \\
A_{ml}^T & 0 & 0 & 0 & 0 & 0 & 0 & 0 & 0 & 0 \\
0 & 0 & 0 & 0 & 0 & 0 & 0 & I_m & 0 & 0 \\
0 & 0 & 0 & 0 & 0 & I_{ml}  & 0 & 0 & 0 & 0 \\
0 & 0 & 0 & 0 & 0 & 0 & 0 & 0 & I_{hm} & 0 \\
0 & 0 & 0 & 0 & 0 & 0 & 0 & 0 & 0 & I_{hw}  \\
A_w^T & 0 & 0 & 0 & 0 & 0 & 0 & 0 & 0 & 0 \\
A_{mc}^T & 0 & 0 & 0 & 0 & 0 & 0 & 0 & 0 & 0 \\
A_{hm}^T & 0 & 0 & 0 & 0 & 0 & 0 & 0 & 0 & 0 \\
A_{hw}^T & 0 & 0 & 0 & 0 & 0 & 0 & 0 & 0 & 0 \\
\end{array}\right)
\end{eqnarray*}
and
\begin{eqnarray*}
F_{21} = 
\left(\begin{array}{cccccccccccc} 
0 & 0 & 0 & 0 & 0 & 0  & 0 & 0 & K_0 & 0& 0 & 0 \\
I_c & 0 & 0 & 0 & 0 & 0 & 0 & 0 & 0 & 0 & 0 & 0 \\
0 & I_{mc} & 0 & 0 & 0 & 0 & 0 & 0 & 0 & K_1 & 0 & 0 \\
0 & 0 & 0 & 0 & 0 & 0 & 0 & 0 & 0 & 0 & 0 & 0 \\
0 & 0 & I_l & 0 & 0 & 0 & 0 & 0 & 0 & 0 & 0 & 0 \\
0 & 0 & 0 & I_{ml} & 0 & K_2 & 0 & 0 & 0 & 0 & 0 & 0 \\
0 & 0 & 0 & 0 & 0 & 0 & 0 & 0 & 0 & 0 & 0 & 0 \\
0 & 0 & 0 & 0 & K_3 & 0 & 0 & 0 & 0 & 0 & 0 & 0 \\
0 & 0 & 0 & 0 & 0 & 0 & K_4 & 0 & 0 & 0 & K_5 & 0 \\
0 & 0 & 0 & 0 & 0 & 0 & 0 & K_6 & 0 & 0 & 0 & K_7 \\
\end{array}\right) 
\end{eqnarray*}
with $K_0 = A_w \partial \zeta / \partial \varphi_w$,
$K_1 = -\partial \omega / \partial \varphi_{mc}$,
$K_2 = -\partial \theta / \partial q_{ml}$,
$K_3 = \partial \eta / \partial q_m$,
$K_4 = \partial \psi / \partial q_{hm}$,
$K_5 = \partial \psi / \partial \varphi_{hm}$,
$K_6 = -\partial \xi / \partial q_{hw}$,
$K_7 = -\partial \xi / \partial \varphi_{hw}$. 

Now, the projector $Q_1$ may be chosen to have the structure
\begin{eqnarray*}
Q_1 & = &
\left(\begin{array}{cc} 
0 & Q_a \\
0 & Q_b
\end{array}\right),
\end{eqnarray*}
where $Q_a$ and $Q_b$ are
\begin{eqnarray*}
%Q_a=
\left(\begin{array}{cccccccccc} 
0 & \hat{Q}_{11} & \hat{Q}_{12}  & \hat{Q}_{13}  & 0 & 0 & 0 & 0 & 0 & 0 \\
0 & \hat{Q}_{21} & \hat{Q}_{22}  & \hat{Q}_{23}  & 0 & 0 & 0 & 0 & 0 & 0 \\
A_l^T \bar{Q} & 0 & 0 & 0 & 0 & 0 & 0 & 0 & 0 & 0 \\
A_{ml}^T \bar{Q}  & 0 & 0 & 0 & 0 & 0 & 0 & 0 & 0 & 0 \\
0 & 0 & 0 & 0 & 0 & 0 & 0 & 0 & 0 & 0 \\
0 & 0 & 0 & 0 & 0 & 0 & 0 & 0 & 0 & 0 \\
0 & 0 & 0 & 0 & 0 & 0 & 0 & 0 & 0 & 0 \\
0 & 0 & 0 & 0 & 0 & 0 & 0 & 0 & 0 & 0 \\
0 & 0 & 0 & 0 & 0 & 0 & 0 & 0 & 0 & 0 \\
0 & 0 & 0 & 0 & 0 & 0 & 0 & 0 & 0 & 0 \\
0 & 0 & 0 & 0 & 0 & 0 & 0 & 0 & 0 & 0 \\
0 & 0 & 0 & 0 & 0 & 0 & 0 & 0 & 0 & 0 
\end{array}\right), \
%\end{eqnarray*}
%and
%\begin{eqnarray*}
%Q_b=
\left(\begin{array}{cccccccccc} 
\bar{Q} & 0 & 0 & 0 & 0 & 0 & 0 & 0 & 0 & 0 \\
0 & \hat{Q}_{11} & \hat{Q}_{12}  & \hat{Q}_{13}  & 0 & 0 & 0 & 0 & 0 & 0 \\
0 & \hat{Q}_{21} & \hat{Q}_{22}  & \hat{Q}_{23}  & 0 & 0 & 0 & 0 & 0 & 0 \\
0 & \hat{Q}_{31} & \hat{Q}_{32}  & \hat{Q}_{33}  & 0 & 0 & 0 & 0 & 0 & 0 \\
0 & 0 & 0 & 0 & 0 & 0 & 0 & 0 & 0 & 0 \\
0 & 0 & 0 & 0 & 0 & 0 & 0 & 0 & 0 & 0 \\
0 & 0 & 0 & 0 & 0 & 0 & 0 & 0 & 0 & 0 \\
0 & 0 & 0 & 0 & 0 & 0 & 0 & 0 & 0 & 0 \\
0 & 0 & 0 & 0 & 0 & 0 & 0 & 0 & 0 & 0 \\
0 & 0 & 0 & 0 & 0 & 0 & 0 & 0 & 0 & 0 
\end{array}\right).
\end{eqnarray*}
Here $\bar{Q}$ is a projector onto $\ke (A_c \ A_{mc} \ A_u \ A_g \
A_w \ A_r \ A_m \ A_{hm} \ A_{hw})^T$, being non-trivial in the
presence of cutsets defined by inductors, meminductors and/or current sources,
whereas
\begin{eqnarray*}  %\label{hatQ}
\hat{Q} = \left( \begin{array}{ccc}
\hat{Q}_{11} & \hat{Q}_{12} & \hat{Q}_{13} \\
\hat{Q}_{21} & \hat{Q}_{22} & \hat{Q}_{23} \\
\hat{Q}_{31} & \hat{Q}_{32} & \hat{Q}_{33}
\end{array}
\right)
\end{eqnarray*}
is a projector onto $\ke (A_c \ A_{mc} \ A_u)$, which does not vanish in
the presence of loops composed of capacitors, memcapacitors and/or
voltage sources.

Additionally, the matrix $F_1=F(I-Q)$ has the expression
\begin{eqnarray*}
\left(\begin{array}{cc} 
0 & 0 \\
F_{21} & 0
\end{array}\right).
\end{eqnarray*}
This gives $E_2 = E_1 - F_1 Q_1$ the form
\begin{eqnarray*}
E_2 & = &
\left(\begin{array}{cc} 
I & -F_{12} \\
0 & -F_{22}  -F_{21}Q_a
\end{array}\right)
\end{eqnarray*}
and, except for the $-$ sign, the lower-right block (which characterizes
the non-singularity of $E_2$) reads
\begin{eqnarray*}
F_{22} + F_{21}Q_a = 
\left(\begin{array}{cccccccccc} 
A_g G A_g^T + A_w W A_w^T & A_c & A_{mc} & A_u & A_l & A_{ml} & A_r & A_m & A_{hm} & A_{hw} \\
-C A_c^T & \hat{Q}_{11} & \hat{Q}_{12}  & \hat{Q}_{13} & 0 & 0 & 0 & 0 & 0 & 0 \\
-C_{\mathrm m} A_{mc}^T & \hat{Q}_{21} & \hat{Q}_{22}  & \hat{Q}_{23} & 0 & 0 & 0 & 0 & 0 & 0 \\
-A_u^T & 0 & 0 & 0 & 0 & 0 & 0 & 0 & 0 & 0 \\
A_l^T \bar{Q} & 0 & 0 & 0 & -L & 0 & 0 & 0 & 0 & 0 \\
A_{ml}^T \bar{Q} & 0 & 0 & 0 & 0 & -L_{\mathrm m} & 0 & 0 & 0 & 0  \\
-A_r^T & 0 & 0 & 0 & 0 & 0 & R & 0 & 0 & 0 \\
-A_m^T & 0 & 0 & 0 & 0 & 0 & 0 & M & 0 & 0 \\
-A_{hm}^T & 0 & 0 & 0 & 0 & 0 & 0 & 0 & M_{\mathrm h} & 0 \\
-W_{\mathrm h}A_{hw}^T & 0 & 0 & 0 & 0 & 0 & 0 & 0 & 0 & I_{hw}
\end{array}\right).
\end{eqnarray*}
Now, by means of a Schur reduction the non-singularity of this matrix is easily 
proved equivalent to that of
\begin{eqnarray*}
%F_{22} + F_{21}Q_a = 
\left(\begin{array}{ccc} 
A_1 M_1 A_1^T + A_2 M_2 A_2^T \bar{Q}& A_3 & A_4 \\ 
-M_3 A_3^T & \tilde{Q} & \check{Q} \\
-A_4^T & 0 & 0 
\end{array}\right)
\end{eqnarray*}
with $A_1 = (A_g \ A_w \ A_r \ A_m \ A_{hm} \ A_{hw})$, $M_1=$ 
block-diag$(G, W, R^{-1}, M^{-1}, M_{\mathrm h}^{-1}, W_{\mathrm h})$,
$A_2 = (A_l \ A_{ml})$, $M_2 =$ block-diag$(L^{-1}, L_{\mathrm m}^{-1})$,
$A_3 = (A_c \ A_{mc})$, $M_3 =$ block-diag$(C, C_{\mathrm m})$, $A_4 = A_u$, 
and
\begin{eqnarray*}
\tilde{Q} =
\left(\begin{array}{cc}
\hat{Q}_{11} & \hat{Q}_{12}  \\
\hat{Q}_{21} & \hat{Q}_{22}  
\end{array}\right), \ 
\check{Q} 
=
\left(\begin{array}{c}
\hat{Q}_{13} \\
\hat{Q}_{23} 
\end{array}\right).
\end{eqnarray*}
The form of this matrix amounts to that arising in the index-two analysis of 
Augmented Nodal Analysis models of nonlinear circuits without memristive
devices, which is proved in \cite[Theorem 5.1(2)]{wsbook}
to be non-singular provided that $M_1$, $M_2$, $M_3$ are positive definite.
In this case, the definiteness of these matrices follows from the 
assumption that the circuit matrices are positive definite and the
proof is complete.
%(WSBOOK p 223, eq (5.50)), 

\hfill $\Box$

\

\noindent In problems without memristive devices, 
Theorems \ref{th-i1} and \ref{th-i2} particularize to the results obtained in 
\cite{et00, wsbook, carentop}.
Our results also extend the characterization
derived in \cite{memristors} for cases including Chua's
%charge-controlled or flux-controlled 
memristors. The scope of the
general index characterization here presented
covers all types of 
devices with differential order one, and therefore applies also to circuits
with $q$- and $\varphi$-memristors,
voltage-controlled memcapacitors, current-controlled meminductors, and
hybrid memristors. These results should be useful in future analytical or
numerical studies of general first order circuits.

\section{Concluding remarks}
\label{sec-con}

We have presented in this paper a comprehensive taxonomy of
a variety of devices which have arisen in nonlinear circuit theory in the
last few years. This taxonomy is organized around the notions of 
the {\em differential} and the {\em state} order of a
device. An exhaustive list of possibly nonlinear
devices with differential order one has been
discussed: besides capacitors and inductors,
these include $q$- and $\varphi$-memristors,
memcapacitors, meminductors, and also the
{\em hybrid memristors} here introduced. 
Hybrid memristors  display a characteristic relating 
all four fundamental circuit variables, and account for devices in
which memory effects of resistive, capacitive and inductive nature coexist.
All these devices are discussed using fully nonlinear characteristics, and
particularize to Chua's memristors and to the memcapacitors and
meminductors of Di Ventra {\em et al.}, respectively, 
when the corresponding 
constitutive relations are linear in certain
variables. A detailed analysis of the differential-algebraic index of
circuits including all possible types of first order devices is also included.

Many aspects remain open and define lines for future investigation.
These include numerical issues, modelling aspects
involving e.g.\ branch-oriented systems and hybrid analysis, or dynamical
properties related to the nature of these circuits' 
operating points, their stability, bifurcations, %phenomena, 
as well as the eventual existence and characterization of 
periodic solutions, oscillations or chaotic 
effects.
Higher order devices and mem-systems are also in the scope of future research.

\section*{Acknowledgements}

The author gratefully acknowledges several suggestions and
bibliographic remarks from Professors M.\ Di Ventra 
(University of California, San Diego) and C.\ Tischendorf 
(University of Cologne, Germany).
%Drive the author's attention to ...

\end{document}